\documentclass[11pt,twoside,reqno]{amsart}
\usepackage{epsfig}

\calclayout

\newtheorem{theorem}[subsection]{\bf Theorem}

\newtheorem{proposition}[subsection]{\bf Proposition}
\newtheorem{cor}[subsection]{\bf Corollary}

\newtheorem{remark}[subsection]{\bf Remark}

\begin{document}

\vspace{1 cm}

\title{An Upper Bound for the regularity of ideals of Borel type.}

\author{Sarfraz Ahmad$^*$, Imran Anwar$^*$}

\thanks{ {The authors are highly grateful to the School of Mathematical Sciences,
  GC University, Lahore, Pakistan in supporting and facilitating this research.
  Our thanks goes also to Prof. Dorin Popescu for his encouragement and valuable suggestions.
  \\$^*$School of Mathematical Sciences, 68-B New Muslim Town, Lahore,
         PAKISTAN.\\\,\,email:\,\,\,sarfraz11@gmail.com, iimrananwar@gmail.com}}

\begin{abstract}
We show that the regularity of monomial ideals of $K[x_1,\ldots,x_n]$ \linebreak ($K$ being a field), whose associated prime ideals are totally ordered by inclusion is upper bounded by a linear function in $n$.

\vskip 0.4 true cm
 \noindent
\end{abstract}
\maketitle

\begin{center}
\textbf{Introduction}
\end{center}

Let $K$ be an infinite field, $S=K[x_1,...,x_n],n\geq 2$ the
polynomial ring over $K$ and $I\subset S$ a monomial ideal. Let
$G(I)$ be the minimal set of monomial generators of $I$ and $deg(I)$
the highest degree of a monomial of $G(I)$. Given a monomial $u\in
S$ set $m(u)=max\{i\mid x_i|u\}$ and $m(I)=max_{u\in G(I)}\ m(u)$.
If $\beta_{ij}(I)$ are the graded Betti numbers of $I$
  then the  {\it regularity} of $I$
is given by $reg(I)=max\{j-i|\beta_{ij}(I)\neq 0\}.$

  Bayer and Mumford \cite{BM}, Caviglia and Sbarra
\cite{CS} and Mayr and Meyer \cite{MaMe} showed that the regularity
of a homogeneous ideal could grow exponentially with respect to its
degree. Set $q(I) = m(I)(deg(I)-1)+1$. Note that if $I$ is a
monomial ideal of height $n$ then $m(I)=n$. In \cite[Remark 2.5]{CS}
it is shown that $reg(I)\leq q(I)$ for a homogeneous ideal $I$ of
height $n$. A monomial ideal $I$ is {\it stable} if for each
monomial $u\in I$ and $1\leq j<m(u)$ it follows $x_ju/x_{m(u)}\in
I$. If $I$ is stable then $reg(I)=deg(I)$ as Eliahou-Kervaire proved
in \cite{EK}. If $I$ is a $p$-Borel ideal then $reg(I)\leq n deg(I)$
as Popescu proved in \cite{Po}.

Let $I$ be a monomial ideal whose associated prime ideals are
totally ordered by inclusion. The purpose of our note is to show
that the regularity of these ideals is bounded by $q(I)$, Corollary
\ref{main}. As can be easily seen, by an appropriate change of
variables, such ideals are the monomial ideals of Borel type
introduced by Herzog, Popescu and Vladoiu in \cite{HPV}. We also
give a new characterization for when an ideal is of Borel type,
Theorem \ref{char} . We  want to point out that Cimpoeas in a recent
preprint \cite{Ci} using different methods achieves the same bound
for a special class of Borel type ideals. We should mention that our
presentation was improved by the kind suggestions of the Referee.

\section{Stable properties of some monomial ideals.}

Let $K$ be a field and  $S=K[x_1,x_2,...,x_n]$. Let $I\subset S$ be
a monomial ideal and $I_{\geq q(I)}$ be the ideal generated by the
monomials of $I$ of degree $\geq q(I)$.

\begin{proposition}
If $I,J$ are monomial ideals such that $I_{\geq q(I)}$ and $J_{\geq q(J)}$ are stable ideals, then
$(I\cap J)_{\geq q(I\cap J)}$ is stable.
\end{proposition}

\begin{proof}
Note that $deg(I\cap J) \geq max\{deg(I),deg(J)\}$,  $m(I\cap J)
\geq max\{m(I),m(J)\}$ and so $q(I\cap J) \geq max\{q(I),q(J)\}$.
Given $v\in (I\cap J)_{\geq q(I\cap J)}$ we get $\frac {x_j\cdot
v}{x_{m(v)}}\in I_{\geq q(I)}\cap J_{\geq q(J)}$ by the stable
property, which is enough.
\end{proof}

\begin{proposition}
A monomial ideal $I=(x_{i_1}^{a_1},...,x_{i_r}^{a_r})\subset S$,
$1\leq i_1<...<i_r\leq n$, $a\in \mathbb{N}^r$ has $I_{\geq q(I)}$ stable if
and only if $i_j=j$ for all $j=1,...,r$.
\end{proposition}

\begin{proof}
$"\Rightarrow"$ Suppose that there exists $1\leq k <i_r$ such that
$k\not\in \{i_1,...,i_r\}$. We have
$x_{i_r}^{q(I)}=x_{i_r}^{q(I)-a_r}x_{i_r}^{a_r}\in I$. As $I_{\geq
q(I)}$ is stable we see by recurrence that $\frac
{x_k^e}{x_{i_r}^e}x_{i_r}^{q(I)}$ for all $1\leq e \leq q(I)$. Thus
$x_k^{q(I)}\in I_{\geq q(I)}$ which is a contradiction since
$I_{\geq q(I)}=\sum_{j=1}^r(x_1,...,x_n)^{q(I)-a_j}\,x_{i_j}^{a_j} $
contains
no pure powers of the variable $x_k$.\\

$"\Leftarrow"$ We show that if  $I=(x_1^{a_1},\ldots,x_r^{a_r})$,
where $a_i\in \mathbb{N}^r$, then $I_{\geq q(I)}$ is stable. Let
$u\in I_{\geq q(I)}$. From above we get $u=v\cdot x_j^{a_j}$ for
some $1\leq j\leq r$ and \linebreak $v\in
(x_1,\ldots,x_n)^{q(I)-a_j}$. If $m(u)>j$ then $\frac
{x_k}{x_{m(u)}}.u=(\frac {x_k\cdot v}{x_{m(u)}})\cdot x_j^{a_j}\in
I_{\geq q(I)}$ for all $k<m(u)$.

If $m(u)=j$ then $u$ belongs to the stable ideal
$(x_1,\ldots,x_r)^{q(I)}$ and it is enough to show that
$(x_1,\ldots,x_r)^{q(I)}\subseteq I_{\geq q(I)}$. Let $w\in
(x_1,\ldots,x_r)^{q(I)}$ then  $w\in I_{q(I)}$, i.e. $w=v.x_i^{a_i}$
where $v\in (x_1,\ldots,x_r)^{q(I) -a_i}$. Indeed, since
$w=x_1^{\alpha_1}x_2^{\alpha_2}\cdots x_r^{\alpha_r}$ with all
$\alpha_i\geq 0$ and $\Sigma_{i=1}^{r} \alpha_i = q(I)$, there
exists some $k$ such that $1\leq k\leq r$ with $\alpha_k\geq a_k$.
So we can take above
$v=x_1^{\alpha_1}...x_k^{\alpha_k-a_k}...x_r^{\alpha_r}$.
\end{proof}

\begin{remark}
In general one cannot get $I_{\geq q(I)-1}$ stable when
$I=(x_1^{a_1},x_2^{a_2},...,x_r^{a_r})$. For example, if $n=2$ and
$I=(x_1^2,x_2^2)$ then $q(I)=3$ and clearly $I_{\geq 2}$ is not stable.
\end{remark}
\vspace{1 pt}

\section{Ideals of Borel Type}
\vspace{10 pt}

A monomial ideal $I$ is {\it of Borel type} if
$(I:x_j^\infty)=(I:(x_1,x_2,...,x_j)^\infty)$ for all $j=1,2,...n.$
These ideals were defined  in \cite{HPV}), where the following was
shown:

\begin{proposition} Let $I$ be a monomial ideal. Then the following conditions
are equivalent:
\begin{enumerate}
    \item $I$ is of Borel type.
    \item If $u\in I$ is a monomial and $1\leq \textrm{i}\leq n $ such that
$x_\textrm{i}^q/u$ for some $q>0$ then for all $1\leq j<\textrm{i}$
there exists an integer $t$ such that $x_j^t(u/x_\textrm{i}^q)\in I$.
\end{enumerate}
\end{proposition}

Our main result is the following:

\begin{theorem}\label{char}
Let $I\in S$ be a monomial ideal.Then the following statements are
equivalent:
\begin{enumerate}
    \item $I$ is an ideal of Borel type.
    \item Each $p\in Ass(S/I)$ has the form $p=(x_1,x_2,...,x_r)$ for some $1\leq r\leq n$.
    \item $I_{\geq q(I)}$ is stable.
\end{enumerate}
\end{theorem}

\begin{proof}

The equivalence of (1) and (2) is well known and is given for
example in the Proposition 5.2 of \cite{HP}

 $(2) \Rightarrow (3)$
Let $I=\bigcap_{e=1}^s J_e$ be the unique irredundant decomposition
of $I$ such that $J_e$ are irreducible monomial ideals (see
\cite[Theorem 5.1.17]{Vi}). Then $J_e=
(x_1^{a_1},\ldots,x_r^{a_r})$,$a_i\in \mathbb{N}^r$ so by
Proposition 1.2,  $(J_e)_{\geq q(J_e)}$ is stable. Thus $(3)$ holds
by Proposition $1.1$.

$(3) \Rightarrow (1)$. Let $u$ be a monomial of $I$ and $1\leq  i\leq n$ such that
$x_i^q|u$ for some $q>0$. By Proposition $2.1$ it is enough to show
that given $1\leq j< i$ there exists an integer $t$ such that
$x_j^t(u/x_i^q)\in I$. We may suppose that $deg(u)\leq q(I)$. Let
$u=x_1^{\alpha_1}\cdots x_n^{\alpha_n}$ and set
$v=x_1^{\alpha_1}\cdots x_{i-1}^{\alpha_{i-1}}x_i^{\alpha_i-q}$,
$w=x_{i+1}^{\alpha_{i+1}}\cdots x_n^{\alpha_n}$, $p=deg(w)+q$. Since
$I_{\geq q(I)}$ is stable and $x_j^{q(I)-deg(u)}u\in I_{\geq q(I)}$,
we get $\frac {x_j}{x_{m(u)}}(x_j^{q(I)-deg(u)}u)\in I_{\geq q(I)}$
and by recurrence we obtain
$x_j^{q(I)-deg(u)+p}\cdot v =\frac{x_j^p}{wx_i^q}(x_j^{q(I)-deg(u)}u)\in
I_{\geq q(I)}$. Thus for $t=q(I)-deg(u)+p$ we have $x_j^t\frac
{u}{x_i^q}=x_j^tvw\in I_{\geq q(I)}.$
\end{proof}

Next we recall a proposition from \cite{ERT}.

\begin{proposition}
\label{ert}
Let $I$ be a monomial ideal and $e\geq deg(I)$ an integer such that
$I_{\geq e}$ is stable. Then $reg(I)\leq e$.
\end{proposition}

We immediately get the following corollaries.
\begin{cor}
If $I$ is of Borel type then $reg(I)\leq m(I)(deg(I)-1)+1$.
\end{cor}

\begin{proof}By the previous theorem we have $I_{\geq q(I)}$
stable, and as $q(I)\geq deg(I)$ we get $reg(I)\leq q(I)$ by
Proposition \ref{ert}.
\end{proof}

\begin{cor}\label{main} If $I$ is a monomial ideal with $Ass(S/I)$ totally
ordered by inclusion then $reg(I)\leq m(I)(deg(I)-1)+1$.
\end{cor}


\medskip
\begin{thebibliography}{0}

\bigskip
\medskip


\bibitem[BM]{BM} D. Bayer, D. Mumford. What can be computed in
Algebraic Geometry? in "{\it Computational Algebraic Geometry and
Commutative Algebra}", Symposia Mathematica XXXIV (1993), 1-48.


\bibitem[CS]{CS} G. Caviglia, E. Sbarra, "{\it Characteristic-free bounds for
the Castelnuovo Mumford regularity}", Compositio Math. 141(2005),
no. 6, 1365-1373.


\bibitem[Ci]{Ci} M. Cimpoeas, "{\it Monomial Ideals with Linear Upper Bound Regularity}",
Preprint, 2006, Arxiv:math.AC/0611064.

\bibitem[EK]{EK} S. Eliahou, M. Kervaire. "{\it Minimal resolutions of some
monomial ideals}" J. Algebra 129 (1990), 1-25.

\bibitem[ERT]{ERT} D. Eisenbud, A. Reeves, B. Totaro, "{\it Initial ideals of
Veronese Subrings}", Adv. in Math., 109(1995), 168-187.

\bibitem[HP]{HP} J. Herzog, D. Popescu, "{\it Finite filtrations of modules
and shellable multicomplexes}", Manuscripta Math., 121(2006),
385-410.

\bibitem[HPV]{HPV}  J. Herzog, D. Popescu, M. Vladoiu, "{\it On the Ext-modules
of ideals of Borel type}", Contemporary Math. 331 (2003), 171-186.

\bibitem[MaMe]{MaMe} E. Mayr, A. Meyer. "{\it The complexity of the word
problem for commutative semigroups and polynomial ideals}", Adv. in
Math 46 (1982), 305-329.

\bibitem[Po]{Po} D. Popescu, "{\it Extremal Betti Numbers and regularity of Borel type
ideals}", Bull. Math. Soc. Sc. Math. Roumanie, Volume 48(96) No. 1,
2005, 65-72.

\bibitem[Vi]{Vi} R. H. Villarreal, "{\it Monomial Algebra}", Marcel Dekker,
Inc, New York, 2001.
\end{thebibliography}
\end{document}